\documentclass[twocolumn]{autart}    % Enable this line and disable the 
\usepackage{graphicx}          % Include this line if your 
\usepackage{algorithm,algorithmic,amsmath,amsfonts,verbatim,mathtools,enumerate,circuitikz,breakcites,todonotes,accents,enumitem,subcaption,balance}
\newcommand{\diag}{\mathop{\bf diag}}

\newcommand{\argmin}{\mathop{\rm argmin}}

\newcommand{\norm}[1]{\left\lVert#1\right\rVert}
\newcommand{\mnorm}[1]{{\left\vert\kern-0.25ex\left\vert\kern-0.25ex\left\vert #1 
    \right\vert\kern-0.25ex\right\vert\kern-0.25ex\right\vert}}
    
\newtheorem{theorem}{Theorem}

\newtheorem{lemma}{Lemma}
\newtheorem{corollary}{Corollary}[theorem]
\newtheorem{remark}{Remark}
\newtheorem{assumption}{Assumption}

\newcommand{\ie}{{\it i.e.}}

\allowdisplaybreaks
\begin{document}

\begin{frontmatter}

\title{Mass-spring-damper Networks for\\ Distributed Optimization in Non-Euclidean Spaces}%\thanksref{footnoteinfo}} % Title, preferably not more than 10 words.
% 
%\thanks[footnoteinfo]{This paper was not presented at any IFAC 
%meeting.}

\author{Yue Yu*},
\author{Beh\c{c}et A\c{c}\i kme\c{s}e},
\author{Mehran Mesbahi}
\address{Department of Aeronautics and Astronautics, University of Washington, Seattle, WA, 98195
}

\thanks{The authors are with the Department of Aeronautics \& Astronautics in the University of Washington, Seattle, WA. (email:{\tt\{yueyu, behcet, mesbahi\}@uw.edu}). This paper was not presented at any IFAC
meeting. Corresponding author Y.~Yu. }
% Please supply full addresses here.

%%The research of M. Mesbahi has been supported by AFOSR grant FA9550-16-1-0022 and DARPA Lagrange grant FA8650-18-2-7836.
          
                             % keyword list or with the 
                                          % help of the Automatica 
                                          % keyword wizard

\begin{abstract}
We consider the problem of minimizing the sum of non-smooth convex functions in non-Euclidean spaces, e.g., probability simplex, via only local computation and communication on an undirected graph. We propose two algorithms motivated by mass-spring-damper dynamics. The first algorithm uses an explicit update  that computes subgradients and Bregman projections only, and matches the convergence behavior of centralized mirror descent. The second algorithm uses an implicit update that solves a penalized subproblem at each iteration, and achieves the iteration complexity of \(\mathcal{O}(1/K)\). The results are also demonstrated via numerical examples.
\end{abstract}

\end{frontmatter}

\section{Introduction}
Given a connected graph, distributed optimization aims to optimize the sum of locally accessible cost functions via local computation and communication \cite{bertsekas1989parallel,boyd2011distributed}. Distributed optimization has a variety of engineering applications, such as formation control \cite{mesbahi2010graph}, distributed tracking and localization \cite{li2002detection}, distributed estimation \cite{accikmecse2014decentralized,lesser2012distributed} and averaging \cite{xiao2007distributed}. Due to such wide range of applications, distributed optimization has been an active area of research during the past two decades. Numerous distributed optimization algorithms have been developed, both in continuous time ~\cite{wang2010control,wang2011control,gharesifard2014distributed,kia2015distributed,qiu2016distributed,zeng2017distributed,yang2017multi,hatanaka2018passivity} and discrete time domains \cite{nedic2009distributed,nedic2010constrained,boyd2011distributed,wei2012distributed,meng2015proximal}. The common feature of these algorithms revolves around sutiable generalization of centralized optimization algorithms to distributed scenarios.

Recently, there has been a growing interest in distributed optimization in non-Euclidean spaces, where the design variable is typically a probability distribution~\cite{dekel2012optimal,levine2016end,gholami2016decentralized,yahya2017collective}.  In order to effectively exploit the structure of such non-Euclidean geometries, several attempts have been made to generalize distributed optimization algorithms from Euclidean to non-Euclidean setting. In particular, the distributed mirror descent method \cite{raginsky2012continuous,li2016distributed,yuan2018optimal,wang2018distributed,doan2019convergence} generalizes the distributed subgradient method \cite{nedic2009distributed}; the distributed dual averaging algorithm \cite{duchi2012dual} generalizes projected distributed subgradient method \cite{nedic2010constrained}; the Bregman parallel direction method of multipliers (BPDMM) \cite{yu2018bregman} generalizes the proximal distributed alternating direction method of multipliers (ADMM) \cite{meng2015proximal}. Compared with their counterparts in Euclidean cases \cite{nedic2009distributed,nedic2010constrained,meng2015proximal}, the key feature of these algorithms is to use a Bregman divergence instead of a quadratic function as the distance generating function, leading to an improved 
complexity bound of \(\mathcal{O}(n/\ln n)\) \cite{wang2014bregman,yu2018bregman} where $n$ represents the dimension of the problem instance.

However, there are still open questions along this line of research. In continuous time domain, compared with distributed Euclidean case \cite{wang2010control,wang2011control,gharesifard2014distributed,kia2015distributed,qiu2016distributed,zeng2017distributed,yang2017multi} and centralized non-Euclidean case \cite{krichene2015accelerated,wibisono2016variational}, the ordinary differential equations (ODE) for distributed non-Euclidean optimization have attracted much less attention. In particular, it is unclear, to the best of our knowledge, whether the ODE setup in \cite{wang2010control} and \cite{raginsky2012continuous} can be jointly utilized to design novel distributed optimization algorithms. In discrete time domain, although BPDMM \cite{yu2018bregman} provides an extension to the proximal distributed ADMM \cite{meng2015proximal}, it requires a computationally expensive mirror averaging step. Moreover, like other algorithms based on ADMM \cite{wei2012distributed,meng2015proximal,yu2018bregman}, BPDMM uses implicit discretization, which requires solving an optimization problem at each iteration. As such, it remains unclear why explicit discretization, which only computes subgradients \cite{duchi2012dual,li2016distributed,doan2019convergence}, cannot achieve similar convergence properties.

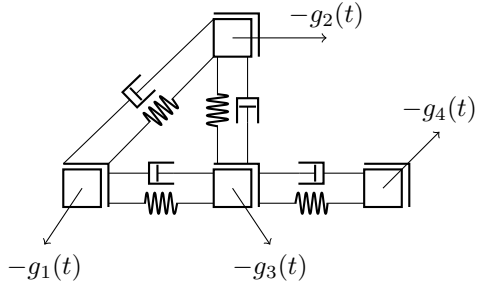
\begin{figure}[ht!]
\centering
\ctikzset{bipoles/length=0.8cm}
\begin{tikzpicture}[scale=0.5]
   \draw[thick] (0,0)rectangle (1, 1);
   \draw[thick] (0, 1.15) -- (1.2, 1.15) --(1.2, 0);
   \draw[->] (0.5, 0.5) to (1.5, -1) node[below] {$-g_3(t)$};
   \draw (1.2, 0.1) to[spring] (4, 0.1);
   \draw (1.2, 0.9) to[damper] (4, 0.9);
   \draw[thick] (4, 0)rectangle (5, 1);
   \draw[thick] (4, 1.15) -- (5.2, 1.15) --(5.2, 0);
   \draw[->] (4.5, 0.5) to (6, 2) node[above] {$-g_4(t)$};
   \draw (0.1, 1.15) to[spring] (0.1, 4);
   \draw (0.9, 1.15) to[damper] (0.9, 4);
   \draw[thick] (0, 4)rectangle (1, 5);
   \draw[thick] (0, 5.15) -- (1.2, 5.15) --(1.2, 4);
   \draw[->] (0.5, 4.5) to (3.0, 4.5) node[above] {$-g_2(t)$};
   \draw (0, 0.1) to[spring] (-2.8, 0.1);
   \draw (0, 0.9) to[damper] (-2.8, 0.9);
   \draw[thick] (-4, 0)rectangle (-3, 1);
   \draw[thick] (-4, 1.15) -- (-2.8, 1.15) -- (-2.8, 0);
   \draw[->] (-3.5, 0.5) to (-4.5, -1) node[below] {$-g_1(t)$};
   \draw (0, 4) to[spring] (-2.8, 1.15);
   \draw (0, 5) to[damper] (-4, 1.15);
   \end{tikzpicture}
\caption{An example of mass-spring-damper network.}
\label{fig: msd}
\end{figure}

Motivated by these questions, as well as the connections between algorithm design and physics \cite{alvarez2000minimizing,alvarez2002second,wibisono2016variational}, we propose two novel algorithms for non-Euclidean distributed optimization with non-smooth convex objective functions over undirected graphs. Our algorithm is based on a mass-spring-damper network model (see Figure~\ref{fig: msd} for an illustration). In particular, we first propose a novel continuous time mass-spring-damper dynamics based on Bregman divergence type kinetic energy for distributed optimization. Using an explicit discretization, we show that discrete time mass-spring-damper dynamics matches the convergence behavior of centralized mirror descent. Further, if an implicit dicretization is used, such convergence can be improved to achieve \(\mathcal{O}(1/K)\) iteration complexity. Finally, we demonstrate our results with numerical experiments.    

Our results extend the existing literature as follows: 1) our continuous time model generalizes the second order ODE model of Euclidean cases \cite{wang2010control} and the first order mirror descent ODE model of non-Euclidean cases \cite{raginsky2012continuous} by combining their attractive features; 2) our discrete time algorithms not only provide a novel extension to distributed mirror descent, but also generalize the proximal function used in distributed ADMM \cite{meng2015proximal} from quadratic functions to a Bregman divergence. Such a generalization is empirically demonstrated  to have faster convergence  than subgradient based algorithms \cite{duchi2012dual,li2016distributed,doan2019convergence}; it also achieves the same iteration complexity and convergence properties as BPDMM \cite{yu2018bregman} with a much more efficient implementation. 

The rest of the paper is organized as follows. \S\ref{sec: preliminaries} provides necessary backgrounds in graph theory and convex analysis. \S\ref{sec: method} first introduces continuous time mass-spring-damper dynamics for distributed optimization, then establishes its convergence in discrete time using different discretization schemes. \S\ref{sec: experiments} compares our algorithms against existing methods via numerical examples. In~\S\ref{sec: conclusion} concluding remarks and comments on future research directions are provided.

\section{Notation and Preliminaries}\label{sec: preliminaries}

Let \(\mathbb{N}\) denote the positive integers, \(\mathbb{R}\) (\(\mathbb{R}_+\)) denote the (non-negative) real numbers, and  \(\mathbb{R}^n\) denote the \(n\)-dimensional real Euclidean space. We use \(\cdot^\top\) to designate matrix (vector) transpose; \(\langle u, v\rangle= u^\top v\) denotes the inner product of two vectors, and \(\norm{\cdot}_2\) is the \(l_2\) norm, \ie , \(\norm{v}_2=\sqrt{\langle v, v\rangle}\) for \(v\in\mathbb{R}^n\). Let \(\diag(v)\in\mathbb{R}^{n\times n}\) denote the diagonal matrix whose diagonal elements are given by vector \(v\in\mathbb{R}^n\); \(\otimes\) designates the Kronecker product. Lastly, \(I_n\in\mathbb{R}^{n\times n}\) denotes the identity matrix and \(\mathbf{1}_n\in\mathbb{R}^n\) the vector of all \(1\)'s.

\subsection{Graph theory}
An undirected graph \(\mathcal{G}=(\mathcal{V}, \mathcal{E})\) consists of a node set \(\mathcal{V}\) and an edge set \(\mathcal{E}\), where an edge is a pair of distinct nodes in \(\mathcal{V}\). The number of nodes and edges in the graph are denoted by \(|\mathcal{V}|\) and \(|\mathcal{E}|\), respectively. Denote by \(\{i j\}\) an edge between nodes \(i\) and \(j\). For an arbitrary orientation on \(\mathcal{G}\), \ie , each edge has a head and a tail, the \(|\mathcal{V}|\times |\mathcal{E}|\) incidence matrix is denoted by \(E(\mathcal{G})\). The columns of \(E(\mathcal{G})\) are indexed by the edge set \(\mathcal{E}\), and the entry on their \(i\)-th row takes the value ``\(1\)" if node \(i\) is the head of the corresponding edge and ``\(-1\)" if it is its tail, and zero otherwise.
When the graph is connected, the nullspace of \(E(\mathcal{G})^\top\) is spanned by \(\mathbf{1}_{|\mathcal{V}|}\) \cite[Theorem 2.8]{mesbahi2010graph}.

\subsection{Convex analysis}
We say that \(g\) is a subgradient of the convex function \(f:\mathbb{R}^n\to\mathbb{R}\) at \(x\) if \cite[Thm 23.2]{rockafellar2015convex},
\begin{equation}
f(x')-f(x)\geq \langle x'-x, g\rangle, \enskip\forall x'.
\label{def: subgradient}
\end{equation}
We denote by \(\partial f(x)\) the set of subgradients of \(f\) at \(x\). Let \(\delta_\mathcal{X}\) be the indicator function of a closed convex set \(\mathcal{X}\), \ie , \(\delta_\mathcal{X}(x)=0\) if \(x\in\mathcal{X}\) and \(+\infty\) otherwise. Then
\begin{equation}
    \partial \delta_\mathcal{X}(x)=N_{\mathcal{X}}(x)=\{g| \langle g, x'-x\rangle\leq 0, \forall x'\in\mathcal{X}\}.
    \label{def: normal cone}
\end{equation}
If \(\psi\) is convex and continuously differentiable, then the Bregman divergence from \(x'\) to \(x\) associated with \(\psi\) is given by \cite{censor1981iterative}
\begin{equation}
B_\psi(x', x)=\psi(x')-\psi(x)-\langle \nabla \psi(x), x'-x\rangle.
    \label{def: Bregman}
\end{equation}
Using \eqref{def: Bregman} one can verify that for any \(x, x^+, x'\),
 \begin{equation}
     \label{eqn: three point}
     \begin{aligned}
     &B_\psi(x', x^+)-B_\psi(x', x)-B_\psi(x^+, x)\\
     =&\langle \nabla\psi(x^+)-\nabla\psi(x), x'-x^+\rangle.
     \end{aligned}
 \end{equation}
Due to the convexity of \(\psi\), \(B_{\psi}(x', x)\)
is always non-negative. If \(\psi\) is \(1\)-strongly convex, \ie, \(\psi-\frac{1}{2}\norm{\cdot}_2^2\) is convex, then using \eqref{def: subgradient} we can show that \(B_{\psi}(x', x)\geq \frac{1}{2}\norm{x'-x}_2^2\).

\section{Mass-spring-damper network dynamics for distributed optimization}\label{sec: method}
In this section, we propose a mass-spring-damper dynamics \eqref{ode: MSD} for the distributed optimization problem \eqref{opt: dist}, leading to algorithms \eqref{alg: explicit} and \eqref{alg: implicit} using explicit and, respectively, implicit discretization. We then establish their convergence properties in Theorem~\ref{thm: explicit} and Theorem~\ref{thm: implicit}.

We consider the distributed optimization problem over the graph \(\mathcal{G}=(\mathcal{V}, \mathcal{E})\) of the following form,
\begin{equation}
\begin{array}{ll}
\underset{x}{\mbox{minimize}} & f(x)=\sum\limits_{i\in\mathcal{V}} f_i(x_i)\\
\mbox{subject to} & E_{s}(\mathcal{G})^\top x=0, \enskip x\in\mathcal{X}=\mathcal{X}_0^{|\mathcal{V}|},
\end{array}
\label{opt: dist}\tag{P}
\end{equation}
where \(x=[x_1^\top, \ldots, x_{|\mathcal{V}|}^\top]^\top\), \(x_i\in\mathbb{R}^n\) for all \(i\in\mathcal{V}\);  \(\mathcal{X}_0^{|\mathcal{V}|}\) is the Cartesian product of \(|\mathcal{V}|\) copies of a closed convex set \(\mathcal{X}_0\subset\mathbb{R}^n\); \(f_i:\mathcal{X}_0\to\mathbb{R}\) is a cost function available to node \(i\) only; matrix \(E_s(\mathcal{G})\) is defined as
\begin{equation}
   E_s(\mathcal{G})=(E(\mathcal{G})\diag(\sqrt{s}))\otimes I_n,
    \label{eqn: incidence}
\end{equation}
where \(\sqrt{s}\) is the element-wise square root of \(s\in\mathbb{R}^{|\mathcal{E}|}\).

To solve problem \eqref{opt: dist}, we propose the following mass-spring-damper dynamics,
\begin{equation}
\begin{aligned}
    \frac{d}{dt}\nabla\psi(x(t))& =-L_d(\mathcal{G})x(t)-E_s(\mathcal{G})u(t)-g(t)\\
   \frac{d}{dt}u(t)&= E_s(\mathcal{G})^\top x(t),
    \end{aligned}\label{ode: MSD}\tag{MSD}
\end{equation}
where \(\psi(x)=\sum_{i\in\mathcal{V}}\psi_0(x_i)\) with \(\psi_0\) being a continuously differentiable convex function over \(\mathcal{X}_0\), \(g(t)=[g_1(t)^\top, \ldots, g_{|\mathcal{V}|}^\top]^\top\) with \(g_i(t)\in\partial f_i(x_i(t))\) for all \(i\in\mathcal{V}\), and finally the weighted Laplacian is given by 
\begin{equation}    \label{eqn: Laplacian}
L_d(\mathcal{G})=(E(\mathcal{G})\diag(d)E(\mathcal{G})^\top)\otimes I_n.
\end{equation}
Dynamics in \eqref{ode: MSD} can be interpreted as a mass-spring-dynamic network as follows. Each node \(i\) in \(\mathcal{V}\) denotes a mass with velocity \(x_i(t)\) experiencing external force \(g_i(t)\); each edge denotes a spring with a spring constant \(s_{\{ij\}}\) and damper with a damping constant \(d_{\{ij\}}\), together connecting node \(i\) and \(j\); see Figure~\ref{fig: msd} for an illustration. If \(\psi_0=\frac{m}{2}\norm{\cdot}_2^2\), then \eqref{ode: MSD} describes the Newton's second law of such mass-spring-damper network, which is also analyzed in \cite{wang2010control,gharesifard2014distributed,hatanaka2018passivity}. Here we use a more general choice of \(\psi_0\) to capture non-quadratic kinetic energy function used in Bregman Lagrangian dynamics \cite{wibisono2016variational}. 

We now group our assumptions on \eqref{opt: dist} and \eqref{ode: MSD} under Assumption~\ref{asp: basic}.
\begin{assumption}\begin{enumerate}
    \item \(\mathcal{G}=(\mathcal{V}, \mathcal{E})\) is undirected and connected. Edge weights \(s, d\in\mathbb{R}^{|\mathcal{E}|}\) are element-wise strictly positive. Let \(\frac{1}{\gamma}=\max\{\alpha|\, \alpha s -d\leq 0\}\) and \(\lambda\) be the largest eigenvalue of \(L_d(\mathcal{G})\).
    \item \(\mathcal{X}_0\subset\mathbb{R}^n\) is a closed convex set, \(\psi_0:\mathcal{X}_0\to\mathbb{R}\) is continuously differentiable and \(1\)-strongly convex, \ie, \(\psi_0-\frac{1}{2}\norm{\cdot}_2^2\) is convex.
    \item \(f_i:\mathcal{X}_0\to\mathbb{R}\) is proper closed convex function for all \(i\in\mathcal{V}\). There exists \(x^\star, u^\star, g^\star\in\partial f(x^\star)\) such that
    \begin{subequations}
    \begin{align}
    E_s(\mathcal{G})^\top x^\star&=0, \enskip x^\star\in\mathcal{X},\label{kkt: primal}\\
    -g^\star-E_s(\mathcal{G})u^\star&\in N_{\mathcal{X}}(x^\star).\label{kkt dual}
    \end{align}
\end{subequations}
\end{enumerate}\label{asp: basic}
\end{assumption}

\begin{remark}
Examples of \(\mathcal{X}_0\) and \(\psi_0\) that satisfy Assumption~\ref{asp: basic} include: 1) \(\mathcal{X}_0=\mathbb{R}^n\) and \(\psi_0(y)=\frac{1}{2}\norm{y}_2^2\); 2) \(\mathcal{X}_0\) is the probability simplex and \(\psi_0(y)=\sum_{i=1}^ny[i]\ln y[i]\) for all \(y\in\mathcal{X}_0\), known as the negative entropy function.\footnote{In this case, \(B_{\psi_0}(x', x)\geq \frac{1}{2}\norm{x'-x}^2_1\)~\cite{beck2003mirror}. Since \(\norm{x}_1\geq \norm{x}_2\) for all \(x\in\mathbb{R}^n\), this further implies that \(\psi_0\) is \(1\)-strongly convex.} See \cite[Sec. 4.3]{bubeck2015convex} for further discussion.
\end{remark}

Note that conditions in \eqref{kkt: primal} and \eqref{kkt dual} ensure primal feasibility and, respectively, stationary condition of \eqref{opt: dist}. 

We first make the following observation on \eqref{ode: MSD}.

\begin{lemma}\label{lem: continuous time}
Suppose that Assumption~\ref{asp: basic} holds and there exists a trajectory \((x(t), u(t))\) for \(t\in[0, \infty)\) generated by \eqref{ode: MSD}. Then 
\begin{equation}
   \ell(\overline{x}(T))-\ell(x^\star)+\frac{1}{2}\langle L_d(\mathcal{G})\overline{x}(T), \overline{x}(T)\rangle\leq \frac{V(x(0), u(0))}{T},
    \label{eqn: duality gap}
\end{equation}
for all \(T>0\), where \(\overline{x}(T)=\frac{1}{T}\int_0^Tx(t)dt\), and
\begin{subequations}\label{eqn: Lyapunov & Lagrangian}
\begin{align}
    V(x, u)=&B_\psi(x^\star, x)+\frac{1}{2}\norm{u-u^\star}_2^2,\\
    \ell(x)=&f(x)+\langle E_s(\mathcal{G}) u^\star, x\rangle.
\end{align}
\end{subequations}
\end{lemma}

\begin{pf*}{Proof}
Using the chain rule we can show that,
\begin{equation}
\begin{aligned}
     \frac{d}{dt}V(x, u)
    =& -\langle \frac{d}{dt}\nabla\psi(x), x^\star-x\rangle+\langle\frac{d}{dt}u, u-u^\star\rangle\\
    =&\langle L_d(\mathcal{G})x+E_s(\mathcal{G})u^\star+g, x^\star-x\rangle\\
    \leq &\textstyle \ell(x^\star)-\ell(x)-\frac{1}{2}\langle L_d(\mathcal{G})x, x\rangle,
    \end{aligned}
\end{equation}
for all \(t\), where the last step is due to \eqref{def: subgradient} and the assumption that \(g(t)\in\partial f(x(t))\). Integrating on both sides from \(t=0\) to \(t=T\) we have 
\begin{equation*}
\begin{aligned}
    &\int_0^t \ell(x(t))-\ell(x^\star)+\frac{1}{2}\langle L_d(\mathcal{G})x(t), x(t)\rangle dt\\
    \leq &V(x(0), u(0))-V(x(T), u(T)).
\end{aligned}    
\end{equation*}
Notice that, by construction, \(V(x(T), u(T))\) is non-negative and \(\ell(x)+\frac{1}{2}\langle L_d(\mathcal{G})x, x\rangle\) is convex. Therefore we can drop the last term on the right hand side. The rest of the proof is an application of Jensen's inequality \(\ell(\overline{x}(T))\leq \int_0^T\ell(x(\tau))d\tau\).\qed
\end{pf*}
Lemma~\ref{lem: continuous time} provides a \(\mathcal{O}(1/T)\) bound on the \emph{running duality gap}~\cite{he20121,meng2015proximal} and the disagreement on edges of the graph along any solution trajectory of \eqref{ode: MSD}. Although such a bound hinges on the existence of such a trajectory, requiring a separate existence proof in general,\footnote{When \(\psi=\frac{1}{2}\norm{\cdot}_2^2\), \eqref{ode: MSD} becomes a differential inclusion with maximal monotone maps that admits a continuous solution trajectory \cite[p.~147]{aubin2012differential}.} it does suggest that a similar property might hold for a discretization of \eqref{ode: MSD}, where the existence of a trajectory is easier to establish. In the following, we aim to show the validity of this intuition.

Consider the following discretization of \eqref{ode: MSD},
\begin{equation}
    \begin{aligned}
     \nabla\psi(x^{k+1})-\nabla\psi(x^k)
     = &-\alpha^k(L_{d}(\mathcal{G})x^k +E_{s}(\mathcal{G})u^k+g^{k})\\
     u^{k+1}-u^k= &\alpha^kE_{s}(\mathcal{G})^\top x^{k+1},
    \end{aligned}
    \label{de: MSD}
\end{equation}
which is an explicit discretization combined with a Gauss-Seidel pass from \(x\)-update to \(u\)-update.

Using \cite[Thm 23.5]{rockafellar2015convex}, we can rewrite the \(x\)-update in \eqref{de: MSD} as,
\begin{equation}
\begin{aligned}
x^{k+1}=&\underset{x\in\mathcal{X}}{\argmin}\, \langle -\nabla\psi(x^k)+\alpha^k(w^k+g^{k}), x\rangle+\psi(x)\\
=&\underset{x\in\mathcal{X}}{\argmin}\,\alpha^k\langle g^{k}+w^k, x\rangle+B_\psi(x, x^k),
\end{aligned}
\end{equation}
where \(w^k=L_d(\mathcal{G})x^k+E_d(\mathcal{G})u^k\).  
Therefore, a trajectory that satisfies \eqref{de: MSD} for all \(k\in\mathbb{N}\) can be computed iteratively as,
\begin{equation}
    \begin{aligned}
        w^k=&L_{d}(\mathcal{G})x^k+E_{s}(\mathcal{G})u^k\\
        x^{k+1}=&\underset{x\in\mathcal{X}}{\argmin}\,\, \alpha^k\langle g^k+ w^k, x\rangle+B_\psi(x, x^k)\\
        u^{k+1}=&u^k+\alpha^k E_{s}(\mathcal{G})^\top x^{k+1},
    \end{aligned}\label{alg: explicit}\tag{MSD-ex}
\end{equation}
where \(g^{k}\in\partial f(x^{k})\), \(\psi(x)=\sum_{i\in\mathcal{V}}\psi_0(x_i)\).

To prove the convergence of \eqref{alg: explicit}, we will use the following property of Bregman divergence.

\begin{lemma}
\label{lem: Pythagorean}
Suppose Assumption~\ref{asp: basic} holds. Given a proper closed convex function \(\phi:\mathcal{X}\to\mathbb{R}\), if \(x^+=\argmin_{y\in\mathcal{X}} \,\,\phi(y)+B_\psi(y, x)\) for some \(x\in\mathcal{X}\), then for any \(x'\in\mathcal{X}\),
\begin{equation*}
    \begin{aligned}
    &B_{\psi}(x', x^+)-B_{\psi}(x', x)\\
    \leq &-\frac{1}{2}\norm{x^+-x}_2^2+\phi(x')-\phi(x^+).
    \end{aligned}
\end{equation*}
\end{lemma}
\begin{pf*}{Proof}
Using Theorem 27.4 in \cite{rockafellar2015convex}, we know that \(x^+=\argmin_{y\in\mathcal{X}} \,\,\phi(y)+B_\psi(y, x)\) if and only if there exists \(g\in\partial \phi(x^+)\) such that
\begin{equation}
    0\leq\langle g+\nabla\psi(x^+)-\nabla\psi(x), x'-x^+\rangle, \label{lem2: eqn1}
\end{equation}
for all \(x'\in\mathcal{X}\).
From \eqref{def: subgradient}, \(\langle g, x'-x^+\rangle\leq \phi(x')-\phi(x^+)\). Substitute this and \eqref{eqn: three point} into \eqref{lem2: eqn1} we have
\[B_\psi(x', x^+)-B_\psi(x', x)\leq-B_\psi(x^+, x)+\phi(x')-\phi(x^+).\]
Finally, from Assumption~\ref{asp: basic} we know that \(\psi-\frac{1}{2}\norm{\cdot}_2^2\) is convex. Using \eqref{def: subgradient} and \eqref{def: Bregman} we can show that \(B_\psi(x^+, x)\geq \frac{1}{2}\norm{x^+-x}_2^2\), which completes the proof. \qed
\end{pf*}
In our subsequent analysis, we will use the following inequality: if \(\lambda\) is the largest eigenvalue of \(L_d(\mathcal{G})\), then
\begin{equation}
\begin{aligned}
    &-\langle L_d(\mathcal{G})x, x^+\rangle\\
    \leq&\frac{\lambda}{2}\norm{x^+-x}_2^2-\frac{1}{2}\langle L_d(\mathcal{G})x, x\rangle-\frac{1}{2}\langle L_d(\mathcal{G})x^+, x^+\rangle,
    \end{aligned}\label{eqn: qudratic bound}
\end{equation}
for all \(x, x^+\); this is due to the fact that \(\frac{1}{2}\langle L_d(\mathcal{G})(x^+-x), x^+-x\rangle\leq \frac{\lambda}{2}\norm{x^+-x}_2^2\). 

With these supporting inequalities, we are now ready to establish the convergence of \eqref{alg: explicit}, which is our first main result.

\begin{theorem}\label{thm: explicit}
Suppose that Assumptions~\ref{asp: basic} holds. If \(0<\alpha^k\leq\min\{\frac{1}{2\lambda}, \frac{1}{\gamma}\}\) for all \(k\in\mathbb{N}\) and \(\norm{g+E_s(\mathcal{G})u^\star}_2\leq G\) for any \(g\in\partial f(x)\) where \(x\in\mathcal{X}\), then along the trajectories generated by algorithm \eqref{alg: explicit},
\begin{equation*}
    \begin{aligned}
    &\ell(\overline{x}^{K})-\ell(x^\star)+\frac{1}{2}\langle L_d(\mathcal{G})\overline{x}^{K}, \overline{x}^{K}\rangle\\
    \leq &\frac{V(x^1, u^1)+G^2\sum_{k=1}^K(\alpha^k)^2}{\sum_{k=1}^K\alpha^k},
    \end{aligned}
\end{equation*}
where \(\overline{x}^{K}=\frac{1}{\sum_{k=1}^K\alpha^k}\sum_{k=1}^K\alpha^kx^k\), \(V(x, u)\) and \(\ell(x)\) are defined as in \eqref{eqn: Lyapunov & Lagrangian}. 
\end{theorem}

\begin{pf*}{Proof}
Applying Lemma~\ref{lem: Pythagorean} to the \(x\)-update in \eqref{alg: explicit} we have,
\begin{equation}\begin{aligned}
    &B_\psi(x^\star, x^{k+1})-B_\psi(x^\star, x^k)\\
    \leq & \textstyle-\frac{1}{2}\norm{x^{k+1}-x^k}_2^2+\alpha^k\langle g^k+E_s(\mathcal{G})u^k, x^\star-x^{k+1}\rangle\\
    &-\alpha^k\langle L_d(\mathcal{G})x^k, x^{k+1}\rangle .
    \end{aligned}\label{eqn: thm1 eqn0}
\end{equation}
In addition, applying \eqref{eqn: three point} to function \(\frac{1}{2}\norm{\cdot}_2^2\) and using the \(u\)-update in \eqref{alg: explicit}, we obtain
\begin{equation}
    \begin{aligned}
    &\textstyle\frac{1}{2}\norm{u^{k+1}-u^\star}^2-\frac{1}{2}\norm{u^k-u^\star}^2\\
    =&\textstyle\frac{(\alpha^k)^2}{2}\langle L_{s}(\mathcal{G}) x^{k+1}, x^{k+1}\rangle\\
    &+\alpha^k\langle u^{k}-u^\star, E_{s}(\mathcal{G})^\top (x^{k+1}-x^\star) \rangle,
    \end{aligned}\label{eqn: thm1 eqn1}
\end{equation}
where we have used \eqref{kkt: primal} and \(L_s(\mathcal{G})=E_s(\mathcal{G})E_s(\mathcal{G})^\top\). Adding \eqref{eqn: thm1 eqn1} to \eqref{eqn: thm1 eqn0} and using \eqref{eqn: qudratic bound}, we can show that,
\begin{equation}
\begin{aligned}
    &V(x^{k+1}, u^{k+1})-V(x^k, u^k)\\
    \leq & \textstyle\alpha^k\langle h^k, x^\star-x^{k+1}\rangle-\frac{\alpha^k}{2}\langle L_d(\mathcal{G})x^{k}, x^{k}\rangle\\
    &\textstyle-\frac{1-\alpha^k\lambda}{2}\norm{x^{k+1}-x^k}_2^2-\frac{\alpha^k}{2}\langle L_c(\mathcal{G})x^{k+1}, x^{k+1}\rangle ,
\end{aligned}\label{eqn: thm1 eqn2}
\end{equation}
where \(c=d-\alpha^k s\) and \(h^k=g^k+E_s(\mathcal{G})u^\star\). Observe that
\begin{equation}\label{eqn: thm1 eqn3}
    \begin{aligned}
        \langle h^k, x^\star-x^k\rangle\leq&\ell(x^\star)-\ell(x^k)\\
        \langle h^k, x^k-x^{k+1}\rangle\leq& \textstyle \frac{\epsilon^k}{2}\norm{h^k}_2^2+\frac{1}{2\epsilon^k}\norm{x^{k+1}-x^k}_2^2
    \end{aligned}
\end{equation}
for any \(\epsilon^k>0\), where the first inequality is due to \eqref{def: subgradient} and \(h^k\in\partial \ell(x^k)\), the second inequality is a completion of square.  Let \(\epsilon^k=\frac{\alpha^k}{1-\alpha^k \lambda}\) and sum up the two inequalities in \eqref{eqn: thm1 eqn3},
\begin{equation}
    \begin{aligned}
    &\langle h^k, x^\star-x^{k+1}\rangle\\
    \leq & \textstyle\ell(x^\star)-\ell(x^k)+\alpha^k G^2+\frac{1-\alpha^k\lambda}{2\alpha^k}\norm{x^{k+1}-x^k}_2^2,
    \end{aligned}\label{eqn: thm1 eqn4}
\end{equation}
where we have used the assumption that \(\norm{h^k}_2\leq G\) and the fact that \(\epsilon^k\leq 2\alpha^k\) when \(0<\alpha^k\leq \frac{1}{2\lambda}\). Substitute \eqref{eqn: thm1 eqn4} into \eqref{eqn: thm1 eqn2} we obtain
\begin{equation*}
    \begin{aligned}
    &\textstyle V(x^{k+1},u^{k+1})-V(x^k, u^k)+\frac{\alpha^k}{2}\langle L_c(\mathcal{G})x^{k+1}, x^{k+1}\rangle\\
    \leq & \textstyle\alpha^k(\ell(x^\star)-\ell(x^k)-\frac{1}{2}\langle L_d(\mathcal{G})x^{k}, x^{k}\rangle)+G^2(\alpha^k)^2.
    \end{aligned}
\end{equation*}
Since \(c=d-\alpha^ks>0\) and \(L_c(\mathcal{G})\) is positive semi-definite when \(\alpha^k\leq \frac{1}{\gamma}\), we can drop the last quadratic term on the right hand side, sum it from \(k=1\) to \(k=K\), and obtain the following
\begin{equation*}
    \begin{aligned}
        &\textstyle\sum_{k=1}^K\alpha^k(\ell(x^k)-\ell(x^\star)+\frac{1}{2}\langle L_d(\mathcal{G})x^k, x^k\rangle)\\
        \leq & \textstyle V(x^1, u^1)-V(x^K, u^K)+G^2\sum_{k=1}^K(\alpha^k)^2.
    \end{aligned}
\end{equation*}
Applying Jensen's inequality to the left hand side and dropping the non-negative \(V(x^K, u^K)\) term on the right hand side, we obtain the desired result.
\qed 
\end{pf*}

\begin{remark}
Theorem~\ref{thm: explicit} matches the convergence results for centralized mirror descent \cite[Thm. 4.1]{beck2003mirror}. One can verify that  
\begin{enumerate}
    \item if sequence \(\{\alpha^k\}\) is square summable but not summable, then Theorem~\ref{thm: explicit} implies an asymptotic convergence.
    \item if \(\alpha^k\) is constant for all \(k\in\mathbb{N}\), then Theorem~\ref{thm: explicit} implies a \(\mathcal{O}(1/K)\) convergence to an error bound of \(\alpha G^2\). Similar results have also been shown in~\cite{wang2018distributed} for the stochastic setting.
\end{enumerate}
\end{remark}

If in addition \(f\) is strongly convex, then Theorem~\ref{thm: explicit} bounds the optimality on nodes and consensus on edges, as the following corollary shows. 

\begin{corollary}\label{cor: explicit}
Under conditions of Theorem~\ref{thm: explicit}, if \(f_i-\frac{\mu_i}{2}\norm{\cdot}_2^2\) is convex with \(\mu_i>0\) for all \(i\in\mathcal{V}\), then
\begin{equation*}
\begin{aligned}
    &\sum_{i\in\mathcal{V}}\frac{\mu_i}{2}\norm{\overline{x}_i^K-x_i^\star}_2^2+\sum_{\{i j\}\in\mathcal{E}}\frac{d_{\{ij\}}}{2}\norm{\overline{x}^K_i-\overline{x}^K_j}_2^2\\
    \leq &\frac{V(x^1, u^1)+G^2\sum_{k=1}^K(\alpha^k)^2}{\sum_{k=1}^K\alpha^k},
\end{aligned}
\end{equation*}
where \(\overline{x}^K_i=\frac{1}{\sum_{k=1}^K\alpha^k}\sum_{k=1}^K\alpha^kx_i^k\) for all \(i\in\mathcal{V}\). 
\end{corollary}

\begin{pf*}{Proof}
For any \(x\in\mathcal{X}\),
\begin{equation*}
    \begin{aligned}
       &\ell(x)-\ell(x^\star)\\
       =&f(x)-f(x^\star)+\langle E_s(\mathcal{G})u^\star+g^\star,  x-x^\star\rangle-\langle g^\star,  x-x^\star\rangle\\
       \geq & \textstyle\sum_{i\in\mathcal{V}}\big(f_i(x_i)-f_i(x_i^\star)-\langle g_i^\star,  x_i-x_i^\star\rangle\big)\\
       \geq& \textstyle\sum_{i\in\mathcal{V}}\frac{\mu_i}{2}\norm{x_i-x_i^\star}_2^2,
    \end{aligned}
\end{equation*}
where the first inequality is due to \eqref{def: normal cone} applied to \eqref{kkt dual} and the fact that \(f(x)=\sum_{i\in\mathcal{V}}f_i(x_i)\); the second inequality is due to \eqref{def: subgradient} applied to \(f_i-\frac{\mu_i}{2}\norm{\cdot}_2^2\) for all \(i\in\mathcal{V}\). Let \(x=\frac{1}{\sum_{k=1}^K\alpha^k}\sum_{k=1}^K\alpha^kx^k\), then the rest of the proof is a direct application of Theorem~\ref{thm: explicit} and the fact that \(\langle L_d(\mathcal{G})\overline{x}^K, \overline{x}^K\rangle=\sum_{\{ij\}\in\mathcal{E}}d_{\{ij\}}\norm{\overline{x}_i^K-\overline{x}_j^K}_2^2\).
\qed
\end{pf*}

Unfortunately, discretization \eqref{alg: explicit} is unable to match the continuous time convergence rate provided in \eqref{eqn: duality gap} of Lemma~\ref{lem: continuous time}. This is mainly due to the non-smoothness of \(f\). However, a matching convergence rate can be achieved by a more careful discretization, as we will show next.   

Consider the following modification of \eqref{alg: explicit}, 
\begin{equation}
    \begin{aligned}
        w^k=&L_{d}(\mathcal{G})x^k+E_{s}(\mathcal{G})u^k\\
        x^{k+1}=&\underset{x\in\mathcal{X}}{\argmin}\,\, \alpha f(x)+\alpha\langle w^k, x\rangle+B_\psi(x, x^k)\\
        u^{k+1}=&u^k+\alpha E_{s}(\mathcal{G})^\top x^{k+1},
    \end{aligned}\label{alg: implicit}\tag{MSD-im}
\end{equation}
where \(\psi(x)=\sum_{i\in\mathcal{V}}\psi_0(x_i)\). Compared with \eqref{alg: explicit}, here we use the function \(f\) in the \(x\)-update, rather than its linear approximation \(\langle g^k, x\rangle\). With this more computationally expensive implicit discretization, \eqref{alg: implicit} is able to achieve a \(\mathcal{O}(1/k)\) iteration complexity that matches the continuous time rate in \eqref{eqn: duality gap}. This is proved by the following theorem, which is our second main result.

\begin{theorem}\label{thm: implicit}
Suppose that Assumption~\ref{asp: basic} holds. If \(\alpha=\min\{\frac{1}{\lambda}, \frac{1}{2\gamma}\}\), then along the trajectories generated by algorithm \eqref{alg: implicit},
\begin{equation*}
\ell(\overline{x}^{K})-\ell(x^\star)+\frac{1}{4}\langle L_d(\mathcal{G})\overline{x}^{K}, \overline{x}^{K}\rangle\leq \frac{V(x^1, u^1)}{\alpha K},
\end{equation*}
where \(\overline{x}^{K}=\frac{1}{K}\sum_{k=1}^Kx^{k+1}\), \(V(x, u)\) and \(\ell(x)\) are defined as in \eqref{eqn: Lyapunov & Lagrangian}. 
\end{theorem}
\begin{pf*}{Proof}
Applying Lemma~\ref{lem: Pythagorean} to the \(x\)-update in \eqref{alg: explicit}, we have
\begin{equation*}\begin{aligned}
    &B_\psi(x^\star, x^{k+1})-B_\psi(x^\star)\\
    \leq & \textstyle -\frac{1}{2}\norm{x^{k+1}-x^k}_2^2+\alpha\langle E_s(\mathcal{G})u^k, x^\star-x^{k+1}\rangle\\
    &-\alpha\langle L_d(\mathcal{G})x^k, x^{k+1}\rangle +\alpha(f(x^\star)-f(x)).
    \end{aligned}
\end{equation*}
Adding \eqref{eqn: thm1 eqn1} with \(\alpha^k=\alpha\) to the above inequality, and again using \eqref{eqn: qudratic bound}, we can then show that,
\begin{equation}
\begin{aligned}
    &V(x^{k+1}, u^{k+1})-V(x^k, u^k)\\
    \leq & \textstyle\alpha(\ell(x^\star)-\ell(x^{k+1}))-\frac{\alpha}{2}\langle L_c(\mathcal{G})x^{k+1}, x^{k+1}\rangle\\
    &\textstyle-\frac{1-\alpha\lambda}{2}\norm{x^{k+1}-x^k}_2^2-\frac{\alpha}{2}\langle L_d(\mathcal{G})x^{k}, x^{k}\rangle,
\end{aligned}\label{eqn: thm2 eqn1}
\end{equation}
where \(c=d-\alpha s\). Since \(\alpha\leq\min\{ \frac{1}{\lambda}, \frac{1}{2\gamma}\}\), \(1-\alpha\lambda\geq 0\) and \(c\geq \frac{1}{2}d\), and \eqref{eqn: thm2 eqn1} implies that
\begin{equation*}
\begin{aligned}
   &\textstyle V(x^{k+1},u^{k+1})-V(x^k, u^k)+\frac{\alpha}{2}\langle L_d(\mathcal{G})x^k, x^k\rangle\\
   \leq &\textstyle\alpha\big(\ell(x^\star)-\ell(x^{k+1})-\frac{1}{4}\langle L_d(\mathcal{G})x^{k+1}, x^{k+1}\rangle\big),
   \end{aligned}
\end{equation*}
Since \(L_d(\mathcal{G})\) is positive semi-definite, we can drop the last quadratic term on the left hand side, sum it from \(k=1\) to \(k=K\) and obtain
\begin{equation*}
\begin{aligned}
   &\textstyle\sum_{k=1}^K\alpha(\ell(x^{k+1})-\ell(x^\star)+\frac{1}{4}\langle L_d(\mathcal{G})x^{k+1}, x^{k+1}\rangle)\\
   \leq &V(x^1, x^1)-V(x^{K+1}, u^{K+1})
   \end{aligned}
\end{equation*}
Applying Jensen's inequality to the left hand side and dropping the non-negative term \(V(x^{k+1}, u^{K+1})\) on the right hand side of the last inequality, we obtain the desired results.
\qed
\end{pf*}

\begin{remark}
\label{rem: mirror average}
Theorem~\ref{thm: implicit} shows that \eqref{alg: implicit} matches the convergence of centralized Bregman proximal method \cite{censor1992proximal}, and achieves the same iteration complexity as those in \cite{yuan2018optimal} and \cite{yu2018bregman}. However, results in \cite{yuan2018optimal} require \(f\) being strongly convex; this assumption is relaxed in Theorem~\ref{thm: implicit}. On the other hand, the algorithm in \cite{yu2018bregman} requires a computationally expensive mirror average step; such computation is not needed in \eqref{alg: implicit}.
\end{remark}

Using the identical argument as Corollary~\ref{cor: explicit}, we can now prove the following extension to Theorem~\ref{thm: implicit}.

\begin{corollary}
Under the conditions in Theorem~\ref{thm: implicit}, if \(f_i-\frac{\mu_i}{2}\norm{\cdot}_2^2\) is convex with \(\mu_i>0\) for all \(i\in\mathcal{V}\), then
\begin{equation*}
\begin{aligned}
    &\sum_{i\in\mathcal{V}}\frac{\mu_i}{2}\norm{\overline{x}_i^K-x_i^\star}_2^2+\sum_{\{i, j\}\in\mathcal{E}}\frac{d_{\{i,j\}}}{4}\norm{\overline{x}^K_i-\overline{x}^K_j}_2^2\\
    \leq &\frac{V(x^1, u^1)}{\alpha K},
    \end{aligned}
\end{equation*}
where \(\overline{x}^K_i=\frac{1}{K}\sum_{k=1}^K x_i^{k+1}\) for all \(i\in\mathcal{V}\).
\end{corollary}
\section{Numerical examples}
\label{sec: experiments}
In this section, we compare both \eqref{alg: explicit} and \eqref{alg: implicit} against existing algorithms for distributed optimization in non-Euclidean spaces, including distributed projected subgradient algorithm \cite{nedic2010constrained}, distributed dual averaging algorithm \cite{duchi2012dual}, distributed mirror descent \cite{li2016distributed,doan2019convergence} and BPDMM \cite{yu2018bregman}, over numerical examples.

\begin{figure}[ht]
  \centering
    \includegraphics[width=0.9\linewidth]{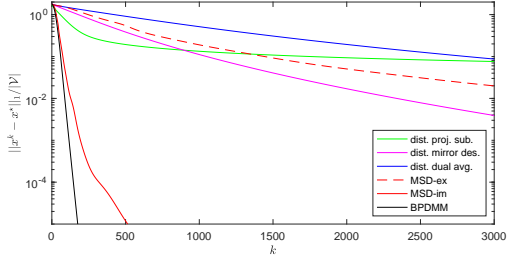}
  \caption{Convergence behavior of various distributed optimization algorithms.} \label{fig: experiment} 
\end{figure}

We choose an instance of problem \eqref{opt: dist} where a) \(\mathcal{G}=(\mathcal{V}, \mathcal{E})\) is randomly generated with \(|\mathcal{V}|=20\) such that each pair of nodes is connected with probability \(0.3\), b) \(f_i(x_i)=\langle a_i, x_i\rangle\) for all \(i\in\mathcal{V}\), where entries of \(a_i\) are sampled uniformly from \([0, 1]\), and c) \(\mathcal{X}_0=\{x\in\mathbb{R}^{10}| x\geq 0, \sum_{i=1}^{10} x[i]=1\}\). For algorithm parameters, we choose \(d=\frac{7}{100}\mathbf{1}_{|\mathcal{E}|}\), and \(s=\lambda d\) where \(\lambda\) is the largest eigenvalue of \(L_d(\mathcal{G})\). For \eqref{alg: explicit}, we let \(\alpha^k=\frac{1}{2\lambda\sqrt{k}}\); for \eqref{alg: implicit}, we let \(\alpha=\frac{1}{2\lambda}\). In comparison, for distributed projected subgradient algorithm, distributed dual averaging, distributed mirror descent,  we use \(P=I-\frac{1}{2+2\Delta} L(\mathcal{G})\) (\(L(\mathcal{G})=E(\mathcal{G})E(\mathcal{G})^\top\), with \(\Delta\) as the largest diagonal element of \(L(\mathcal{G})\) \cite{duchi2012dual}), for the stochastic matrix and \(\alpha^k=\frac{1}{2\lambda\sqrt{k}}\) for the stepsize; for BPDMM, we use \(P=I-\frac{1}{2+2\Delta} L(\mathcal{G})\) for the stochastic matrix and \(\alpha=\frac{1}{2\lambda}\) for the stepsize.\footnote{Here \(\alpha\) corresponds to \(\frac{1}{\rho}\) in \cite{yu2018bregman}.} Except for the distributed projected subgradient, all algorithms use \(\psi_0(y)=\sum_{i=1}^ny[i]\ln y[i]\) for all \(y\in\mathcal{X}_0\). The convergence of these algorithms, using the same randomly generated initialization, are shown in Figure~\ref{fig: experiment}.

Figure~\ref{fig: experiment} shows that \eqref{alg: explicit} shares similar convergence behavior as the distributed mirror descent and the distributed dual averaging. This is mainly due the fact that both algorithms are based on explicit computation of subgradients, which requires diminishing step sizes to ensure convergence. On the other hand, \eqref{alg: implicit} shares similar convergence behavior as BPDMM, as the implicit update allows a constant step size. The price is, in general, a more expensive computation at each iteration. Further, \eqref{alg: implicit} almost achieves the same convergence speed as the BPDMM with a more efficient implementation, as pointed out in Remark~\ref{rem: mirror average}. 
\section{Conclusions}\label{sec: conclusion}

In this paper, we have developed two novel algorithms for distributed optimization in non-Euclidean spaces based on the mass-spring-damper dynamics. Our results not only match those for centralized mirror descent, but also improve previous methods by allowing more efficient implementation. Nevertheless, the proposed algorithm has a few limitations, most notably, their mere applicability to undirected graphs. Our future research direction is  motivated by this limitation, as well as extensions for smooth and stochastic optimization as well as time-varying graphs.

% OR

%\begin{figure}
%\begin{center}
%\epsfig{file=jcaesar,width=7cm}
%\caption{Gaius Julius Caesar, 100--44 B.C.}
%\label{fig1}
%\end{center}
%\end{figure}

%\begin{ack}                               % Place acknowledgements
%Partially supported by the Roman Senate.  % here.
%\end{ack}

\bibliographystyle{apalike}         % Include this if you use bibtex 
%\balance
\bibliography{reference}           % and a bib file to produce the 
                                 % bibliography (preferred). The
                                 % correct style is generated by
                                 % Elsevier at the time of printing.

%\appendix
     % Sections and subsections are supported  
                                        % in the appendices.
\end{document}